\documentstyle[12pt]{amsart}
\def\NZQ{\Bbb}               

\def\ZZ{{\NZQ Z}}
\def\RR{{\NZQ R}}

\def\opn#1#2{\def#1{\operatorname{#2}}} 
\opn\ini{in} \opn\chara{char} \opn\length{\ell}
\opn\projdim{proj\,dim} \opn\injdim{inj\,dim} \opn\rank{rank}
\opn\depth{depth} \opn\grade{grade} \opn\height{height}
\opn\embdim{emb\,dim} \opn\codim{codim}

\opn\Tr{Tr} \opn\bigrank{big\,rank}
\opn\superheight{superheight}\opn\lcm{lcm}
\opn\trdeg{tr\,deg}%
\opn\reg{reg} \opn\link{link} \opn\sta{star}

\opn\div{div} \opn\Div{Div} \opn\cl{cl} \opn\Cl{Cl}
\opn\gr{gr}

\def\qed{\ifhmode\textqed\fi
   \ifmmode\ifinner\quad\qedsymbol\else\dispqed\fi\fi}
\def\textqed{\unskip\nobreak\penalty50
    \hskip2em\hbox{}\nobreak\hfil\qedsymbol
    \parfillskip=0pt \finalhyphendemerits=0}
\def\dispqed{\rlap{\qquad\qedsymbol}}

\newtheorem{de}{Definition}[section]

\newtheorem{pr}[de]{Proposition}
\newtheorem{co}[de]{Corollary}

\newtheorem{te}[de]{Theorem}

\begin{document}

\title{Cohen-Macaulay classes which are not conic}
\author {Cornel Baetica}

\address{University of Bucharest, Faculty of Mathematics,
70109 Bucharest, Romania} \email{baetica@@al.math.unibuc.ro}




\begin{abstract}
We provide two examples which show that for an affine semigroup
ring the set of Cohen-Macaulay classes can be larger than the set
of conic classes.
\end{abstract}

\maketitle


\section*{Introduction}
Let $S\subset\ZZ^n$ be a normal semigroup. It can be described as
the set of lattice points in a finite generated rational cone, i.
e. $S=\{x\in\ZZ^n:\sigma_i(x)\geq 0, i=1,\ldots,s\}$ is the set of
lattice points satisfying a system of homogeneous inequalities
given by linear forms $\sigma_i$ with integral coefficients. For a
field $K$ the semigroup ring $R=K[S]$ is normal. Suppose in
addition that the group gp($S$) generated by $S$ equals $\ZZ^n$,
and the presentation of $S$ is irredundant. Bruns and Gubeladze
\cite{BG} have studied the divisorial ideals of $R$ of the form
$KT$, where $T=\{z\in\ZZ^n:\sigma_i(z)\geq\sigma_i(\beta)\}$, for
some $\beta\in\RR^n$. They called these divisorial ideals {\em
conic}, and it is known that conic divisorial ideals are
Cohen-Macaulay (see Stanley \cite{St} and Dong \cite {D}). The set
of conic classes contains all torsion classes in the divisor class
group, and it is strictly larger whenever $\Cl(R)$ is non-torsion.
In this note we provide two examples which show that the set of
Cohen-Macaulay classes can be larger than the set of conic
classes.

\section{The Segre product of three polynomial rings}
It is known that all Cohen-Macaulay classes are conic when the
ring $R$ is the Segre product of two polynomial rings over a
field, see \cite{BG}. We wonder whether the same property holds
when one considers the Segre product of three polynomial rings. We
can prove that in this case the number of Cohen-Macaulay classes
is greater than the number of conic classes.

Consider a field $K$ and the Segre product $R=K[X_iY_jZ_k: 1\leq
i\leq m, 1\leq j\leq n, 1\leq k\leq p]$ of polynomial rings
$R_1=K[X_1,\ldots,X_m]$, $R_2=K[Y_1,\ldots,Y_n]$ and
$R_3=K[Z_1,\ldots,Z_p]$, $m,n,p\geq 2$, with its standard
embedding in
$$P=K[X_1,\ldots,X_m,Y_1,\ldots,Y_n,Z_1,\ldots,Z_p].$$ We often
denote the Segre product of two rings $A$ and $B$ as $A\#B$. The
ring $R$ can be described as a semigroup ring $K[S]$ with
$S=\{(\alpha,\beta,\gamma)\in\ZZ_+^{m+n+p}:
|\alpha|=|\beta|=|\gamma|\}$, where
$|\alpha|=\alpha_1+\ldots+\alpha_m$ for $\alpha\in\ZZ^m$. It has
the divisor class group isomorphic to $\ZZ^2$ and $P$ decomposes
as an $R$-module into a direct sum of rank $1$ $R$-modules $M_c$,
$c\in\ZZ^2$, such that  $M_c$ is isomorphic to a divisorial ideal
of class $c$, see \cite[Theorem 2.1]{BG}. We can write
$$M_{(i,j)}=\sum_{|\alpha|=|\beta|+i,\;
|\alpha|=|\gamma|+j}KX^{\alpha}Y^{\beta}Z^{\gamma},$$ and the
crucial point is that $M_{(i,j)}$ is also a Segre product of some
shifts of the polynomial rings above, namely
$M_{(i,j)}=R_1\#R_2(-i)\#R_3(-j)$. To decide when the modules
$M_{(i,j)}$ are Cohen-Macaulay we will use a criterion of
St\"uckrad and Vogel \cite[Theorem]{SV} and we will record it
here. However some definitions are needed.
\begin{de}
Let $R$ be a positive graded $K$-algebra, and $M$ a finitely
generated graded $R$-module. We define $a(M)$ as being the degree
of $H_M(t)$, the Hilbert series of $M$ and $r(M)= {\rm
inf}\{n\in\ZZ:M_n\neq 0\}$.
\end{de}
Let us remark that $a(M)$ is an extension of the $a$-invariant for
modules.
\begin{te}\label{cr}{\rm (St\"uckrad and Vogel \cite{SV})}
Let $R_1$ and $R_2$ be positively graded $K$-algebras and $M_1$,
$M_2$ Cohen-Macaulay graded $R_1$-resp. $R_2$-modules of Krull
dimension greater or equal than 2. Then $M_1\#M_2$ is a
Cohen-Macaulay module (over $R_1\#R_2$) if and only if
$a(M_1)+1\leq r(M_2)$ and $a(M_2)+1\leq r(M_1)$.
\end{te}
Another result we will use in the sequel is the following
\begin{pr}{\rm (Bruns and Guerrieri \cite{BGr})}
Let $R_1=K[X_1,\ldots,X_m]$ and $R_2=K[Y_1,\ldots,Y_n]$. Then
$R_1\#R_2(-i)$ is Cohen-Macaulay if and only if $-(m-1)\leq i\leq
n-1$.
\end{pr}
We are ready now to prove the following
\begin{pr}
The ring $R$ has $(m^2+n^2+p^2)+(mn+mp+np)-2(m+n+p)+1$
Cohen-Macaulay classes, and $(mn+mp+np)-(m+n+p)+1$ conic classes.
\end{pr}
\begin{pf} To decide when the modules $M_{(i,j)}$ are Cohen-Macaulay
we have to consider the following cases

Case 1:  $i\geq 0$, $j\geq i$. If $i\leq n-1$, then $R_1\#R_2(-i)$
is Cohen-Macaulay and the module $M_{(i,j)}$ is Cohen-Macaulay if
and only if $a(R_1\#R_2(-i))+1\leq j$ and $a(R_3(-j))+1\leq i$. As
$a(R_1\#R_2(-i))=-{\rm max}(m,n-i)$, then we get $M_{(i,j)}$
Cohen-Macaulay if and only if $j-i\leq p-1$.

If $j\leq p-1$, then $R_1\#R_3(-j)$ is Cohen-Macaulay  and the
module $M_{(i,j)}$ is Cohen-Macaulay if and only if
$a(R_1\#R_3(-j))+1\leq i$ and $a(R_2(-i))+1\leq j$. As
$a(R_1\#R_3(-j))=-{\rm max}(m,p-j)$, then we get $M_{(i,j)}$
Cohen-Macaulay if and only if $j\leq p-1$.

If $i\geq n$ and $j\geq p$, then we consider the following
situations: for $j-i\leq p-1$ we have $R_2\#R_3(i-j))$
Cohen-Macaulay and $M_{(i,j)}$ is isomorphic to
$R_1\#(R_2\#R_3(i-j))(-i)$. Applying again the criterion \ref{cr}
we get an impossible inequality, so $M_{(i,j)}$ can not be
Cohen-Macaulay in this case. If $j-i\geq p$, then $\mu(M_{(i,j)})$
the minimal number of generators of $M_{(i,j)}$ is greater than
$e(R)$ the multiplicity of $R$ and by Serre's numerical
Cohen-Macaulay criterion (see \cite[Theorem 4.7.11]{BH}) we get
again that $M_{(i,j)}$ is not Cohen-Macaulay.

It is not hard to see that the number of Cohen-Macaulay classes is
$pn+(p-n)(p-n+1)/2$ in this case.

Case 2: $i\geq 0$ and $0\leq j<i$. As before, for $i\leq n-1$ we
get $M_{(i,j)}$ Cohen-Macaulay for all $i$, for $j\leq p-1$ the
module $M_{(i,j)}$ is Cohen-Macaulay if and only if $i-j\leq n-1$,
and for $i\geq n$ and $j\geq p$ the module $M_{(i,j)}$ is not
Cohen-Macaulay. We get $p(n-1)$ Cohen-Macaulay classes in this
case.

Case 3: $i\geq 0$ and $j<0$. The same argument as above leads us
to the following: for $-j\leq m-1$ $M_{(i,j)}$ is Cohen-Macaulay
if and only if $i\leq n-1$, $M_{(i,j)}$ is Cohen-Macaulay for
$i-j\leq n-1$, and $M_{(i,j)}$ is not Cohen-Macaulay for $-j\geq
m$ and $i-j\geq n$.  The number of Cohen-Macaulay classes is
$n(m-1)+(n-m)(n-m+1)/2$ in this case.

The following three cases are symmetric to the first three.

Case 4: $i<0$ and $i\geq j$. The module $M_{(i,j)}$ is
Cohen-Macaulay  for any $-j\leq m-1$, for $i-j\leq n-1$ the module
$M_{(i,j)}$ is Cohen-Macaulay if and only if $-i\leq m-1$, and
$M_{(i,j)}$ is not Cohen-Macaulay for $-j\geq m$ and $i-j\geq n$.
Therefore the number of Cohen-Macaulay classes is $n(m-1)$ in this
case.

Case 5: $i<0$ and $i<j\leq 0$. The module $M_{(i,j)}$ is
Cohen-Macaulay  for any $-i\leq m-1$, for $j-i\leq p-1$ the module
$M_{(i,j)}$ is Cohen-Macaulay if and only if $-j\leq m-1$, and
$M_{(i,j)}$ is not Cohen-Macaulay for $-i\geq m$ and $j-i\geq p$.
We get $m(p-1)$ Cohen-Macaulay classes in this case.

Case 6: $i<0$ and $j>0$. The module $M_{(i,j)}$ is Cohen-Macaulay
for any $j-i\leq p-1$, for $-i\leq m-1$ the module $M_{(i,j)}$ is
Cohen-Macaulay if and only if $j\leq p-1$, and $M_{(i,j)}$ is not
Cohen-Macaulay for $-i\geq m$ and $j-i\geq n$. The number of
Cohen-Macaulay classes in this case is
$(m-1)(p-1)+(p-m-1)(p-m)/2$.

Summing up we get $(m^2+n^2+p^2)+(mn+mp+np)-2(m+n+p)+1$
Cohen-Macaulay classes.

On the other hand, it is not hard to count how many conic classes
we have. It relies on a remark made by Bruns and Gubeladze
\cite{BG} in the proof of their Proposition 3.6. In fact, any
conic class corresponds to the class of $\wp_1^{(i)}\cap
\wp_2^{(j)}$, where $\wp_1$ is the divisorial ideal of $R$
generated by the products $X_1Y_kZ_l$, $1\leq k\leq n$, $1\leq
l\leq p$, and $\wp_2$ is the divisorial ideal of $R$ generated by
the products $X_kY_1Z_l$, $1\leq k\leq m$, $1\leq l\leq p$. To get
all conic classes it is enough to consider $i=\lceil{a-b}\rceil$
and $j=\lceil{a-c}\rceil$ with $-m<a\leq 0$, $-(n-1)<b\leq 0$, and
$-(p-1)<c\leq 0$. An easy counting provides the number of conic
classes: $(mn+mp+np)-(m+n+p)+1$.
\end{pf}
Finally we get
\begin{co}
For the Segre product of three polynomial rings over a field the
number of Cohen-Macaulay classes is strictly greater than the
number of conic classes.
\end{co}

\section{The Segre product of two Veronese embeddings}
Our second example involves the Segre product of two Veronese
embeddings. Let $K$ be a field and let $K[X_1,\ldots,X_m]$ be a
polynomial ring over $K$. The $c$-th Veronese subring of
$K[X_1,\ldots,X_m]$ is the $K$-subalgebra of $K[X_1,\ldots,X_m]$
generated by the set of the monomials of degree $c$. We denote it
by $K[X_1,\ldots,X_m]^{(c)}$. Let us consider
$R=R_1^{(c)}\#R_2^{(d)}$, where $R_1=K[X_1,\ldots,X_m]$,
$R_2=K[Y_1,\ldots,Y_n]$, $m,n,c,d\geq 1$, $(c,d)=1$, with its
standard embedding in the polynomial ring
$P=K[X_1,\ldots,X_m,Y_1,\ldots,Y_n]$. The ring $R$ can be
described as a semigroup ring $K[S]$ with
$S=\{(\alpha,\beta)\in\ZZ_+^{m+n}: d|\alpha|=c|\beta|\}$. Its
divisor class group is isomorphic to $\ZZ$ and $P$ decomposes as
an $R$-module into a direct sum of rank $1$ $R$-modules $M_i$,
$i\in\ZZ$, such that $M_i$ is isomorphic to a divisorial ideal of
class $i$, see again \cite[Theorem 2.1]{BG}. We can write
$M_i=\sum_{d|\alpha|-c|\beta|=i}KX^{\alpha}Y^{\beta}$ and the
point is that $M_i$ is also a Segre product of two Veronese
embeddings of some shifts of the polynomial rings above, namely
$M_i=R_1(-vi)^{(c)}\#R_2(-ui)^{(d)}$, where $u,v\geq 1$ such that
$cu-dv=1$.
\begin{pr}
The ring $R$ has $m+n+c+d-3$ conic classes, and at least $dm+cn-1$
Cohen-Macaulay classes.
\end{pr}
\begin{pf}
Every conic class corresponds to the class of $p_1^{(i)}\cap
p_2^{(j)}$, where $p_1$ is the divisorial ideal of $R$ generated
by the monomials $x\in S$ such that $\sigma_1(x)\geq 1$ and $p_2$
is the divisorial ideal of $R$ generated by the monomials $x\in S$
such that $\sigma_2(x)\geq 1$, where the $\sigma_i$'s are the
support forms of $S$. To get all conic classes it is enough to
consider $i=\lceil{ca-b}\rceil$ and $j=\lceil{da-b'}\rceil$ with
$-1<a\leq 0$, $-(m-1)<b\leq 0$, and $-(n-1)<b'\leq 0$. Their
number is $m+n+c+d-3$.

The St\"uckrad and Vogel's criterion gives us that $M_i$ is
Cohen-Macaulay if and only if
$-\lceil{(m-vi)/c}\rceil+1\leq\lceil{ui/d}\rceil$ and
$-\lceil{(n-ui)/d}\rceil+1\leq\lceil{vi/c}\rceil$. The first
inequality holds for $i\geq -dm+1$, while the second holds for
$i\leq cn-1$. But in general, it is possible to get values for
$i$, outside of these ranges, that satisfies both inequalities.
\end{pf}
\begin{co}
The ring $R$ has the same number of Cohen-Macaulay, respectively
conic classes if and only if $c=d=1$, $c=m=1$, or $d=n=1$.
\end{co}
\begin{pf}
For $c=1$ we get exactly $dm+n-1$ Cohen-Macaulay classes, while
for $d=1$ their number is $m+cn-1$. On the other side, when the
number of Cohen-Macaulay classes equals the number of conic
classes, then $m+n+c+d-3$ should be greater or equal than
$dm+cn-1$. This leads to the following inequality
$(d-1)(m-1)+(c-1)(n-1)\leq 0$. It holds if and only if
$(d-1)(m-1)=0$ and $(c-1)(n-1)=0$. However, it is not difficult to
see that for $m=n=1$ and $c>1, d>1$ the number of Cohen-Macaulay
classes is greater than $c+d-1$.
\end{pf}

\section*{Acknowledgement} The author is indebted to Prof.
W. Bruns who invited him at the Osnabr\"uck University in
June-July 2001 and would like to thank him for valuable
suggestions. The author also acknowledges the financial support
given by NATO, grant A/01/41484.


\begin{thebibliography}{99}

\bibitem{BG}
W. Bruns and J. Gubeladze, {\em Divisorial linear algebra of
normal semigroup rings}, preprint.

\bibitem{BGr}
W. Bruns and A. Guerrieri, {\em The Dedekind-Mertens formula and
determinantal rings}, Proc. Amer. Math. Soc. {\bf 127}(1999),
657-663. MR {\bf 99f:}13013

\bibitem{BH}
W. Bruns and J. Herzog, {\em Cohen-Macaulay rings}, Cambridge
University Press, Cambridge, 1993. MR {\bf 95h:}13020

\bibitem{D}
Xun Dong, {\em Alexander duality and a theorem of Danilov and
Stanley}, preprint.

\bibitem{St}
R. Stanley, {\em Combinatorics and invariant theory}. In D. K.
Ray-Chaudury (ed.), {\em Relations between combinatorics and other
parts of mathematics}, Proc. Symp. Pure Math. {\bf 34}(1979),
345-355. MR {\bf 80e:}15020

\bibitem{SV}
J. St\"uckrad and W. Vogel, {\em On Segre products and
applications}, J. Algebra {\bf 54}(1978), 374-389. MR {\bf
80c:}14031

\end{thebibliography}
\end{document}